\documentclass[12pt]{article}
\usepackage{amssymb,amsfonts,amsmath}
\usepackage[english]{babel}
\usepackage{color}

\begin{document}

\begin{center}{\large\bf On the convergence of exotic formal series solutions of an ODE. A proof by the implicit mapping theorem}\end{center}
\bigskip

\centerline{R.\,R.\,Gontsov\footnote{Institute for Information Transmission Problems of RAS,
Bolshoy Karetny per. 19, build.1, Moscow 127051 Russia. Moscow Power Engineering Institute, 
Krasnokazarmennaya 14, Moscow 111250 Russia. E-mail: gontsovrr@gmail.com.}, 
I.\,V.\,Goryuchkina\footnote{Keldysh Institite of Applied Mathematics of RAS,
Miusskaya sq. 4, Moscow 125047 Russia. E-mail: igoryuchkina@gmail.com.}}

\begin{abstract}

We propose a sufficient condition of the convergence of a complex power type formal series of the form
$\varphi=\sum_{k=1}^{\infty}\alpha_k(x^{{\rm i}\gamma})\,x^k$, where $\alpha_k$ are functions meromorphic at 
the origin and $\gamma\in{\mathbb R}\setminus\{0\}$, that satisfies an analytic ordinary differential equation 
(ODE) of a general type. An example of a such type formal solution of the third Painlev\'e equation is presented
and the proposed sufficient condition is applied to check its convergence.

\end{abstract}

\section{Introduction}

We consider a non-linear ODE of order $n$,
\begin{eqnarray}\label{eq1}
F(x,y,\delta y,\ldots,\delta^n y)=0,
\end{eqnarray}
where $\delta=x(d/dx)$ and $F(x,y_0,y_1,\ldots,y_n)$ is a holomorphic function near $0\in{\mathbb C}^{n+2}$.
Suppose that the equation \eqref{eq1} possesses a formal series solution $y=\varphi$ of the form
\begin{equation}\label{eq2}
\varphi=\sum\limits_{k=1}^{\infty}\alpha_k(x^{{\rm i}\gamma})\,x^k, \qquad {\rm i}=\sqrt{-1}, \quad \gamma\in\mathbb{R}\setminus\{0\}, 
\end{equation}
 where $\alpha_k(t)$ are meromorphic functions at the origin:
 $$
\alpha_k(t)=t^{-r_k}\sum\limits_{\ell=0}^{\infty}\alpha_{k\ell}\,t^{\ell}, \qquad \alpha_{k\ell}\in\mathbb{C}, 
\qquad r_k\in\mathbb{Z},
$$
with some common punctured disc $\overline{\cal D}\setminus\{0\}=\{0<|t|\leqslant R\}$ of convergence.

The series \eqref{eq2} will be called {\it exotic}, in the terminology of A.\,D.\,Bruno \cite{Br}. In particular, the 
Painlev\'e III, V, VI equations possess formal solutions of such type \cite{BrGo}, \cite{BrPa}, \cite{Gu}, \cite{Sh} 
(in the last paper \cite{Sh} by Shimomura the {\it convergence} of such formal solutions of the fifth Painlev\'e equation 
is established). Thus, our present goal is to obtain some general condition for the convergence of an exotic formal 
series solution of the equation \eqref{eq1}. In this sense, our work continues a series of articles \cite{GG1}, 
\cite{GG3}, \cite{GG4} where similar questions were studied for {\it generalized} formal power series  
$$
\varphi^{\rm pow}=\sum_{k=0}^{\infty}c_kx^{\lambda_k}, \qquad {\rm Re}\,\lambda_0\leqslant{\rm Re}\,\lambda_1\leqslant
\ldots, \quad {\rm Re}\,\lambda_k\rightarrow+\infty,
$$
and for formal {\it Dulac} series
$$
\varphi^{\rm Dul}=\sum_{k=1}^{\infty}P_k(\ln x)x^k, \qquad P_k\in{\mathbb C}[t].
$$ 
Those articles were inspired by the original paper by B.\,Malgrange \cite{Ma} on the Maillet theorem for classical 
formal power series solutions of a non-linear ODE. 

Note that the series (\ref{eq2}) can be written in the form of a complex power type series,
$$
\varphi=\sum\limits_{k=1}^{\infty}\sum\limits_{\ell=0}^{\infty}\alpha_{k\ell}\,x^{k+{\rm i}\gamma(\ell-r_k)},
$$
however, this differs from a generalized power series $\varphi^{\rm pow}$ since for each $k\geqslant1$, there are 
infinite number of power exponents with the same real part $k$ and, therefore, the set of power exponents of such 
a series cannot be ordered with respect to the real part increasing. 

We propose the following sufficient condition of the convergence of (\ref{eq2}), which is the main result of the 
present paper. 
\medskip

{\bf Theorem 1.} {\it Let $(\ref{eq2})$ be a formal solution of the equation $(\ref{eq1}):$
$$
F(x,\Phi)=0, \qquad \Phi:=(\varphi,\delta\varphi,\dots,\delta^n\varphi),
$$
such that $\displaystyle \frac{\partial F}{\partial y_n}(x,\Phi)\not\equiv 0$. Furthermore, let each exotic formal 
series $\displaystyle \frac{\partial F}{\partial y_j}(x,\Phi)$ be of the form
$$
\frac{\partial F}{\partial y_j}(x,\Phi)=a_j(x^{{\rm i}\gamma})x^m+b_j(x^{{\rm i}\gamma})x^{m+1}+\ldots, \qquad 
j=0,1,\ldots,n,
$$
with the same $m$ for all $j$, where $a_n\not\equiv0$ and
$$
{\rm ord}_0\,a_j\geqslant {\rm ord}_0 \,a_n, \qquad j=0,1,\dots,n.
$$
Then in any open sector $S$ with the vertex at the origin, of opening less than $2\pi$ and of sufficiently 
small radius, such that $S\subset\{\arg x>(-1/\gamma)\ln R\}$ if $\gamma>0$ $($and $S\subset\{\arg x<(-1/\gamma)\ln R\}$ 
if $\gamma<0)$, the series $(\ref{eq2})$ converges uniformly in $S$, thus representing there a branch of a holomorphic 
function.}
\medskip

There is also a qualitative difference between generalized power series solutions and exotic series solutions of 
(\ref{eq1}). Since $\varphi^{\rm pow}$ converges in any sufficiently small sector with the vertex at the origin
and of opening less than $2\pi$, with no restriction on its location (see \cite{GG1}), singular points of a solution
represented by $\varphi^{\rm pow}$ cannot accumulate to the origin along a ray or ray-like curve. However, they
can accumulate along spirals (for example, like in the case of Painlev\'e equations \cite{Bre}, \cite{Sh}). On 
the other hand, singular points of a solution represented by (\ref{eq2}) can accumulate to the origin along a ray 
that is not contained in the domain from Theorem 1. Indeed, if $t=a$, $|a|>R$, is a singular point of some 
$\alpha_k(t)$ then the equation $x^{{\rm i}\gamma}=a$ has an infinite number of solutions located on the ray 
$\arg x=(-1/\gamma)\ln |a|$ and having the origin as a limit point (again, an illustration of this situation is
provided by Painlev\'e equations, see \cite{Gu2}). 

We note that the condition of Theorem 1 is similar to the corresponding sufficient conditions of the convergence of 
generalized formal power series solutions and formal Dulac series solutions of (\ref{eq1}) obtained in \cite{GG1} and 
\cite{GG4} (all they generalize Malgrange's condition \cite{Ma} of the convergence of a classical formal power series 
solution and require the linearization of the equation along a formal solution to be {\it Fuchsian} in some sense). 
We prove Theorem 1 in a series of lemmas adapting the main idea of Malgrange, an application of the implicit mapping 
theorem for Banach spaces, to the situation under consideration.

\section{An ODE in a reduced form}

A starting point in the proof of Theorem 1 is a reduction of the initial equation to a special form, which is similar 
to what Malgrange did studying classical formal power series solutions.
\medskip

{\bf Lemma 1.} {\it Under the assumptions of Theorem 1, there exists $N\in\mathbb{N}$ such that the transformation
\begin{equation}\label{eq3}
 y=\sum\limits_{k=1}^N\alpha_k(x^{{\rm i}\gamma})x^k+x^N u
\end{equation}
reduces the equation \eqref{eq1} to an equation of the form
\begin{equation}\label{eq4}
  \sum_{j=0}^na_j(x^{{\rm i}\gamma})\,(\delta+N)^ju=x\,M(x,x^{{\rm i}\gamma},u,\delta u,\dots,\delta^nu),
\end{equation}
where the function $M(x,t,u_0,\dots,u_n)$ is meromorphic in $t$ and holomorphic in the rest of variables, near 
$0\in{\mathbb C}^{n+3}$.}
\medskip

{\bf Proof.} For any $N\in{\mathbb N}$, let us write the formal solution $\varphi$ in the form
$$
\varphi=\sum\limits_{k=1}^N\alpha_k(x^{{\rm i}\gamma})x^k+x^N\psi=:\varphi_N+x^N\psi.
$$
Respectively,
$$
\Phi=\Phi_N+x^N\Psi, \qquad \Psi=(\psi,(\delta+N)\psi,\ldots,(\delta+N)^n\psi).
$$
Then the Taylor formula implies
\begin{eqnarray}\label{taylor}
0&=&F(x,\Phi_N+x^N\Psi)=F(x,\Phi_N)+x^N\sum_{j=0}^n\frac{\partial F}{\partial y_j}(x,\Phi_N)\psi_j+ \nonumber \\
 & &+\frac{x^{2N}}2\sum_{i,j=0}^n\frac{\partial^2 F}{\partial y_i\partial y_j}(x,\Phi_N)\psi_i\psi_j+\ldots,
\end{eqnarray}
where $\psi_j=(\delta+N)^j\psi$.
\medskip

{DEFINITION.} Let us define the valuation of $\varphi$ as
$$
{\rm val}\,\varphi:=\min\{k\mid \alpha_k\not\equiv0\}.
$$
\medskip

Again, by the Taylor formula, 
$$
\frac{\partial F}{\partial y_j}(x,\Phi)-\frac{\partial F}{\partial y_j}(x,\Phi_N)=x^N\sum_{i=0}^n
\frac{\partial^2 F}{\partial y_i\partial y_j}(x,\Phi_N)\psi_i+\ldots,
$$
furthermore ${\rm val}(\psi_i)\geqslant1$ for all $i$. If we choose $N>m$ (recall that the integer $m\geqslant0$ comes 
from the condition of Theorem 1) then there will hold\footnote{There will be another requirement for choosing $N$, which 
is not used in the proof of Lemma 1 and which we explain later, in Lemma 2.}
$$
\frac{\partial F}{\partial y_j}(x,\Phi_N)=a_j(x^{{\rm i}\gamma})x^m+\tilde b_j(x^{{\rm i}\gamma})x^{m+1}+\dots, 
$$
for each $j=0,1,\ldots,n$, that is, the leading coefficient $a_j$ will be preserved if one substitutes the finite 
sum $\Phi_N$ instead of $\Phi$ in $\frac{\partial F}{\partial y_j}$. Now from the relation (\ref{taylor}) it follows 
that 
$$
{\rm val}\,F(x,\Phi_N)\geqslant m+N+1.
$$
Hence, dividing (\ref{taylor}) on $x^{m+N}$, one obtains 
$$
\sum_{j=0}^n a_j(x^{{\rm i}\gamma})(\delta+N)^j\psi=x\,M(x,x^{{\rm i}\gamma},\psi,\delta\psi,\dots,\delta^n\psi),
$$
where $M(x,t,u_0,\dots,u_n)$ is a function meromorphic in $t$ and holomorphic in the rest of variables, near 
$0\in{\mathbb C}^{n+3}$. This finishes the proof of the lemma. {\hfill $\Box$}
\medskip

{\bf Remark 1.} According to the conditions of Lemma 1, one can assume that $a_n\equiv1$ and the rest of $a_j$'s
to be {\it holomorphic} near the origin.
\medskip

One sees that the reduced equation (\ref{eq4}) possesses an exotic formal series solution
\begin{equation}\label{eq7}
\psi=\sum_{k=1}^{\infty}c_k(x^{{\rm i}\gamma})\,x^k, \qquad 
c_k(t)=\alpha_{k+N}(t)=t^{-\nu_k}\sum_{\ell\geqslant0}c_{k\ell}\,t^{\ell},
\end{equation}
moreover, this is its unique formal solution of such a form.
\medskip

{\bf Lemma 2.} {\it The formal series $(\ref{eq7})$ is a unique exotic series satisfying $(\ref{eq4})$.} 
\medskip

{\bf Proof.} The differential operator $\delta$ acts on the series $\psi$ in the following way:
$$
\delta\psi=\sum_{k=1}^{\infty}\bigl(k+{\rm i}\gamma\delta_t\bigr)c_k(t)|_{t=x^{{\rm i}\gamma}}\,x^k, \quad 
\delta_t=t(d/dt),
$$
whence
$$
(\delta+N)^j\psi=\sum_{k=1}^{\infty}\bigl(k+N+{\rm i}\gamma\delta_t\bigr)^jc_k(t)|_{t=x^{{\rm i}\gamma}}\,x^k.
$$
Hence, each coefficient $c_k$ is a solution of the (non-homogeneous) linear differential equation
\begin{equation}\label{fuchs}
\sum_{j=0}^n a_j(t)\bigl(k+N+{\rm i}\gamma\delta_t\bigr)^jc_k(t)=\tilde c_k(t), \quad k=1,2,\ldots,
\end{equation}
where $\tilde c_k$ is a known function meromorphic at the origin, which is determined by the previous 
$\tilde c_1,\ldots,\tilde c_{k-1}$ (for example, $\tilde c_1(t)=M(0,t,0,\ldots,0)$, {\it etc}.).

By assumptions (see Remark 1), the linear differential operator on the left hand side of (\ref{fuchs}) is Fuchsian 
at the point $t=0$. Its exponents at this point are the roots of the polynomial
$$
P(\lambda)=\sum_{j=0}^n a_j(0)\bigl(k+N+{\rm i}\gamma\lambda\bigr)^j.
$$
Assuming additionally $N$ to be large enough in such a way that $P(\lambda)$ has no integer roots (this is the second 
requirement to $N$) we get that any non-trivial local solution of the corresponding homogeneous equation
$$
\sum_{j=0}^n a_j(t)\bigl(k+N+{\rm i}\gamma\delta_t\bigr)^jy(t)=0
$$
ramifies at the Fuchsian singular point $t=0$. Hence, $c_k(t)$ is a unique solution of (\ref{fuchs}) meromorphic at $t=0$.
{\hfill $\Box$}
\medskip

{\bf Lemma 3.} {\it For the pole order $\nu_k$ of each coefficient $c_k(t)$ of the series $(\ref{eq7})$ at $t=0$, 
the following estimate holds: $\nu_k\leqslant k\mu$, where $\mu$ is the pole order of $M(x,t,u_0,\ldots,u_n)$ at $t=0$.}
\medskip

{\bf Proof.} Since the polynomial $P(\lambda)=\sum_{j=0}^n a_j(0)\bigl(k+N+{\rm i}\gamma\lambda\bigr)^j$
has no integer roots (in particular, $P(-\nu_k)\ne0$) from (\ref{fuchs}) it follows that $\nu_k$ is equal to the 
pole order of $\tilde c_k(t)$ at $t=0$. Hence, as $\tilde c_1(t)=M(0,t,0,\ldots,0)$, one has $\nu_1\leqslant\mu$ and 
further proceeds by induction as follows.

The Taylor terms of the function $M$ that contribute to the term $\tilde c_k(x^{{\rm i}\gamma})x^k$ of the exotic 
series $x\,M(x,x^{{\rm i}\gamma},\psi,\delta\psi,\dots,\delta^n\psi)$, are those of the form
$$
c\,x^r(x^{{\rm i}\gamma s})\psi^{k_0}(\delta\psi)^{k_1}\dots(\delta^n\psi)^{k_n}, \qquad 0\leqslant r\leqslant k-1,
$$
with the additional restriction $k_0+k_1+\ldots+k_n\leqslant k-1-r$. Therefore, denoting
$$
\delta^j\psi=\sum_{k=1}^{\infty}\bigl(k+{\rm i}\gamma\delta_t\bigr)^jc_k(t)|_{t=x^{{\rm i}\gamma}}\,x^k=:
\sum_{k=1}^{\infty}c_k^j(x^{{\rm i}\gamma})\,x^k
$$
one sees that $\tilde c_k(t)$ is a linear combination of functions of the form
$$
t^s\bigl(c_{l_1}^0\dots c_{l_{k_0}}^0\bigr)\bigl(c_{m_1}^1\dots c_{m_{k_1}}^1\bigr)\dots
\bigl(c_{n_1}^n\dots c_{n_{k_n}}^n\bigr),
$$ 
furthermore 
$$
s\geqslant-\mu, \qquad \sum_{i=1}^{k_0}l_i+\sum_{i=1}^{k_1}m_i+\dots+\sum_{i=1}^{k_n}n_i\leqslant k-1.
$$
By a natural inductive assumption, for $l<k$ the pole order of each $c^j_l$ at $t=0$, which is equal to $\nu_l$,
does not exceed $l\mu$, hence
$$
\nu_k\leqslant\mu+\sum_{i=1}^{k_0}l_i\mu+\sum_{i=1}^{k_1}m_i\mu+\dots+\sum_{i=1}^{k_n}n_i\mu\leqslant\mu+(k-1)\mu=k\mu.
$$
{\hfill $\Box$}
\medskip

{\bf Remark 2.} It follows from the inductive determination of $c_k$'s as solutions of linear ODEs with the same 
homogeneous part, the initial requirement for them to be holomorphic in the common punctured disc 
$\overline{\cal D}\setminus\{0\}$ is, in fact, an internal property of these functions.

\section{Banach spaces of exotic series}

To prove Theorem 1, one should prove the convergence of the formal exotic series solution $\psi$ of the reduced equation
(\ref{eq4}). The idea is, using the implicit mapping theorem for Banach spaces, to prove the existence of a convergent 
exotic series solution of (\ref{eq4}) and then, from uniqueness proved in Lemma 2, to deduce the convergence of $\psi$.  

For $\nu\geqslant1$, we consider the space ${\cal O}_{\nu}(\overline{\cal D}\setminus\{0\})$ of functions holomorphic 
in a punctured disc and having a pole of order at most $\nu$ at $t=0$. Expanding 
$f\in{\cal O}_{\nu}(\overline{\cal D}\setminus\{0\})$ into a Laurent series, 
$$
f=\frac1{t^{\nu}}\sum_{k\geqslant 0}f_k\,t^k,
$$
we define its norm $\|f\|$ by
$$
\|f\|=\frac1{r^{\nu}}\sum_{k\geqslant 0}|f_k|\,R^k, \qquad 0<r<R.
$$
\medskip

{\bf Lemma 4.} {\it The space ${\cal O}_{\nu}(\overline{\cal D}\setminus\{0\})$ with the norm $\|\cdot\|$ is a 
Banach space.} 
\medskip

{\bf Proof.} Since $\max\limits_{r\leqslant|t|\leqslant R}|f(t)|\leqslant\|f\|$, any fundamental sequence 
$f^{(m)}\in{\cal O}_{\nu}(\overline{\cal D}\setminus\{0\})$ is also fundamental with respect to the norm of 
uniform convergence and, therefore, converges uniformly to a function $F$ holomorphic in the annulus 
$\{r\leqslant|t|\leqslant R\}$: $F=\sum_{k\in{\mathbb Z}}c_k\,t^k$. 
Due to the Cauchy formula for the coefficients of a Laurent series, 
$$
c_k=\frac1{2\pi{\rm i}}\int\limits_{|t|=\rho}\frac{F(t)}{t^{k+1}}\,dt, \qquad r<\rho<R,
$$
each coefficient $c_k$ of $F$ is equal to the limit of the sequence of the corresponding coefficients 
$c_k^{(m)}$ of $f^{(m)}$, hence $c_k=0$ for $k<-\nu$ and $F\in{\cal O}_{\nu}(\overline{\cal D}\setminus\{0\})$. 
Thus, it remains to prove that $\|f^{(m)}-F\|\rightarrow0$ ($m\rightarrow\infty$). Since for any $\varepsilon>0$ 
there exists $m_0$ such that $\|f^{(m)}-f^{(n)}\|\leqslant\varepsilon$ for $m,n\geqslant m_0$:
$$
\frac1{r^{\nu}}\sum_{k\geqslant 0}|c_{k-\nu}^{(m)}-c_{k-\nu}^{(n)}|\,R^k\leqslant\varepsilon,
$$
one has
$$
\frac1{r^{\nu}}\sum_{k=0}^K|c_{k-\nu}^{(m)}-c_{k-\nu}^{(n)}|\,R^k\leqslant\varepsilon
$$ 
for any $K\geqslant0$. Taking the limit in the above inequality with respect to $n\rightarrow\infty$, 
we obtain
$$
\frac1{r^{\nu}}\sum_{k=0}^K|c_{k-\nu}^{(m)}-c_{k-\nu}|\,R^k\leqslant\varepsilon,
$$ 
which holds for every $K\geqslant0$. Therefore,
$$
\frac1{r^{\nu}}\sum_{k\geqslant0}|c_{k-\nu}^{(m)}-c_{k-\nu}|\,R^k\leqslant\varepsilon
$$
for $m\geqslant m_0$. This means that $\|f^{(m)}-F\|\leqslant\varepsilon$ for $m\geqslant m_0$ and, therefore,
$\|f^{(m)}-F\|\rightarrow0$ ($m\rightarrow\infty$). {\hfill $\Box$}
\medskip

Now keeping in mind Lemma 3 we put $\nu_k=k\mu$ and introduce the spaces
\begin{eqnarray*}
H^j&=&\Bigl\{\theta=\sum_{k=1}^{\infty}p_k(x^{{\rm i}\gamma})\,x^k \mid p_k\in{\cal O}_{\nu_k}
(\overline{\cal D}\setminus\{0\}), \; \|\theta\|_j=\sum_{k=1}^{\infty}\|(k+N+{\rm i}\gamma\delta_t)^jp_k\|
<+\infty\Bigr\}, \\
   & & j=0,1,\ldots,n,
\end{eqnarray*}
of exotic series. These turn out to be Banach spaces which further allows us to apply the implicit mapping theorem
for a properly chosen operator in these very spaces, to prove the existence of a convergent exotic series solution
of (\ref{eq4}).
\medskip

{\bf Lemma 5.} {\it Each $H^j$, $j=0,1,\ldots,n,$ is a Banach space.}
\medskip

{\bf Proof.} Let $\{\theta^{(m)}\}$ be a fundamental sequence in $H^j$, that is, for any $\varepsilon>0$ there exists 
$m_0$ such that $\|\theta^{(m)}-\theta^{(n)}\|_j\leqslant\varepsilon$ for $m,n\geqslant m_0$:
\begin{equation}\label{ineqH}
\sum_{k=1}^{\infty}\|(k+N+{\rm i}\gamma\delta_t)^j(p^{(m)}_k-p^{(n)}_k)\|\leqslant\varepsilon.
\end{equation}
It follows that for each fixed $k$ one has
$$
\|(k+N+{\rm i}\gamma\delta_t)^j(p^{(m)}_k-p^{(n)}_k)\|\leqslant\varepsilon,
$$
that is, $\{(k+N+{\rm i}\gamma\delta_t)^jp^{(m)}_k\}_{m=1}^{\infty}$ is a fundamental sequence in 
${\cal O}_{\nu_k}(\overline{\cal D}\setminus\{0\})$. Hence, by Lemma 4 it converges to some 
$q_k\in{\cal O}_{\nu_k}(\overline{\cal D}\setminus\{0\})$. From (\ref{ineqH}) one deduces that for 
any fixed $K\geqslant1$
$$
\sum_{k=1}^K\|(k+N+{\rm i}\gamma\delta_t)^j(p^{(m)}_k-p^{(n)}_k)\|\leqslant\varepsilon
$$
whenever $m,n\geqslant m_0$, therefore
$$
\sum_{k=1}^K\|(k+N+{\rm i}\gamma\delta_t)^jp^{(m)}_k-q_k\|\leqslant\varepsilon
$$
for $m\geqslant m_0$. Since the last inequality holds for any $K$, one has
$$
\sum_{k=1}^{\infty}\|(k+N+{\rm i}\gamma\delta_t)^jp^{(m)}_k-q_k\|\leqslant\varepsilon
$$
for $m\geqslant m_0$. Expanding $q_k$ into the Laurent series in ${\cal D}\setminus\{0\}$,
$$
q_k(t)=\sum_{\ell\geqslant-\nu_k}q_{k\ell}\,t^{\ell},
$$
one sees that 
$$
\tilde q_k(t)=\sum_{\ell\geqslant-\nu_k}\frac{q_{k\ell}}{(k+N+{\rm i}\ell\gamma)^j}\,
t^{\ell}\in{\cal O}_{\nu_k}(\overline{\cal D}\setminus\{0\})
$$
and $(k+N+{\rm i}\gamma\delta_t)^j\tilde q_k=q_k$, that is, $\theta^{(m)}$ converges to 
$\theta=\sum_{k=1}^{\infty}\tilde q_k(x^{{\rm i}\gamma})\,x^k$ in $H^j$. {\hfill $\Box$}
\medskip

{\bf Lemma 6.} {\it The Banach spaces $H^j$ possess the following properties.
\smallskip

{\rm i)} The following inclusions hold:
$$
H^n\subset\ldots\subset H^1\subset H^0\subset{\cal O}(S_1), 
$$
where $S_1$ is any open sector with the vertex at the origin, of opening less than $2\pi$ and of radius $1$,
such that 
\begin{eqnarray*}
S_1 & \subset & \{(-1/\gamma)\ln R<\arg x<(-1/\gamma)\ln r\} \quad\mbox{ if } \quad\gamma>0, \\
S_1 & \subset & \{(-1/\gamma)\ln r<\arg x<(-1/\gamma)\ln R\} \quad\mbox{ if } \quad\gamma<0.  
\end{eqnarray*}
Furthermore, $\|\theta\|_{j-1}\leqslant\|\theta\|_j$ for any $\theta\in H^j$.
\smallskip

{\rm ii)} For any $\theta_1,\theta_2\in H^0$, there hold $\theta_1\theta_2\in H^0$ and 
$\|\theta_1\theta_2\|_0\leqslant\|\theta_1\|_0\|\theta_2\|_0$.
\smallskip

{\rm iii)} For any $\theta\in H^j$ one has $\delta\theta\in H^{j-1}$ and the operator $\delta: H^j\rightarrow H^{j-1}$ 
is continuous.}
\medskip

{\bf Proof.} i) Let $\theta=\sum_{k=1}^{\infty}p_k(x^{{\rm i}\gamma})\,x^k\in H^j$. Expanding
$$
(k+N+{\rm i}\gamma\delta_t)^jp_k=\frac1{t^{\nu_k}}\sum_{\ell\geqslant0}p_{k\ell}\,t^{\ell}, \qquad
\|(k+N+{\rm i}\gamma\delta_t)^jp_k\|=\frac1{r^{\nu_k}}\sum_{\ell\geqslant0}|p_{k\ell}|\,R^{\ell}, 
$$
one has
$$
(k+N+{\rm i}\gamma\delta_t)^{j-1}p_k=\frac1{t^{\nu_k}}\sum_{\ell\geqslant0}
\frac{p_{k\ell}}{k+N+{\rm i}(\ell-\nu_k)\gamma}\,t^{\ell}, 
$$
whence the estimate $\|(k+N+{\rm i}\gamma\delta_t)^{j-1}p_k\|\leqslant(1/k)\|(k+N+{\rm i}\gamma\delta_t)^jp_k\|$ 
follows. This implies $\|\theta\|_{j-1}\leqslant\|\theta\|_j$ and $H^j\subset H^{j-1}$, $j=1,\ldots,n$.

As for the inclusion $H^0\subset{\cal O}(S_1)$, for any $x\in S_1$ one has
$$
|x^{{\rm i}\gamma}|=e^{-\gamma\arg x}\in(r,R),
$$
therefore 
$$
|p_k(x^{{\rm i}\gamma})\,x^k|\leqslant\max\limits_{r\leqslant|t|\leqslant R}|p_k(t)|\leqslant\|p_k\|,
$$
and the convergence of $\sum_{k=1}^{\infty}\|p_k\|$ implies the uniform convergence of 
$\sum_{k=1}^{\infty}p_k(x^{{\rm i}\gamma})\,x^k$ in $S_1$.
\smallskip

ii) For any functions $f\in{\cal O}_{\nu_k}(\overline{\cal D}\setminus\{0\})$,
$g\in{\cal O}_{\nu_{\ell}}(\overline{\cal D}\setminus\{0\})$ one has $\|fg\|\leqslant\|f\|\,\|g\|$.
Indeed, let the Laurent expansions of $f$ and $g$ in $\overline{\cal D}\setminus\{0\}$ be, respectively,
$$
f=\frac1{t^{\nu_k}}\sum_{s\geqslant0}f_s\,t^s, \qquad g=\frac1{t^{\nu_{\ell}}}\sum_{s\geqslant0}g_s\,t^s.
$$
Then
$$
\|fg\|=\frac1{r^{\nu_k+\nu_{\ell}}}\sum_{s\geqslant0}\Bigl|\sum_{i=0}^sf_ig_{s-i}\Bigr|R^s\leqslant
\frac1{r^{\nu_k+\nu_{\ell}}}\sum_{s\geqslant0}\Bigl(\sum_{i=0}^s|f_i|\,|g_{s-i}|\Bigr)R^s=\|f\|\,\|g\|.
$$
Therefore, for $\theta_1=\sum_{k=1}^{\infty}p_k(x^{{\rm i}\gamma})\,x^k$, 
$\theta_2=\sum_{k=1}^{\infty}q_k(x^{{\rm i}\gamma})\,x^k\in H^0$ there holds 
$$
\theta_1\theta_2=\sum_{k=2}^{\infty}\Bigl(\sum_{\ell=1}^{k-1}p_{\ell}(x^{{\rm i}\gamma})
q_{k-\ell}(x^{{\rm i}\gamma})\Bigr)x^k,
$$ 
where $p_{\ell}\,q_{k-\ell}\in{\cal O}_{\nu_k}(\overline{\cal D}\setminus\{0\})$ for each $\ell$
(since $\nu_l+\nu_{k-l}=\nu_k$). Furthermore
$$
\|\theta_1\theta_2\|_0=\sum_{k=2}^{\infty}\Bigl\|\sum_{\ell=1}^{k-1}p_{\ell}\,q_{k-\ell}\Bigr\|
\leqslant\sum_{k=2}^{\infty}\sum_{\ell=1}^{k-1}\|p_{\ell}\|\,\|q_{k-\ell}\|=\|\theta_1\|_0\|\theta_2\|_0.
$$ 

iii) If $\theta\in H^j$ then $(\delta+N)\theta\in H^{j-1}$ and $\|(\delta+N)\theta\|_{j-1}=\|\theta\|_j$, 
hence 
$$
\delta\theta=(\delta+N)\theta-N\theta\in H^{j-1}\quad\mbox{ and }\quad 
\|\delta\theta\|_{j-1}\leqslant(1+N)\|\theta\|_j,
$$
that is, $\delta: H^j\rightarrow H^{j-1}$ is a continuous operator. {\hfill $\Box$}

\section{Finishing the proof of Theorem 1 by the implicit mapping theorem}

Theorem 1, in view of Lemma 2, will follow from the following lemma.
\medskip

{\bf Lemma 7.} {\it The equation $(\ref{eq4})$ possesses an exotic series solution $\theta_0\in{\cal O}(S)$, where $S$ 
is a sector as in Theorem 1.}
\medskip

{\bf Proof.} Consider the mapping $h:\mathbb C\times H^n\rightarrow H^0$ defined by the formula
$$
h(\lambda,\theta)=\sum_{j=0}^n a_j(x^{{\rm i}\gamma})(\delta+N)^j\theta-\lambda x\,M(\lambda x,x^{{\rm i}\gamma}, \theta,\delta\theta,\dots,\delta^n\theta)
$$
in some neighbourhood of the point $(0,0)\in\mathbb C\times H^n$, $h(0,0)=0$. By Lemma 6, this mapping and its 
derivative $\partial h/\partial\theta$ are continuous at $(0,0)$. The local solvability of the equation 
$h(\lambda,\theta)=0$, $\theta=\theta(\lambda)$, and as a consequense, Lemma 7 will follow from the implicit mapping 
theorem (\cite[Th. 10.2.1]{Die}) if we prove that
$$
\frac{\partial h}{\partial\theta}(0,0)=\sum_{j=0}^n a_j(x^{{\rm i}\gamma})(\delta+N)^j: H^n\rightarrow H^0
$$  
is a linear homeomorphism. To prove this, we should first prove that for any sequence $\{q_k(t)\}$, 
$q_k\in{\cal O}_{\nu_k}(\overline{\cal D}\setminus\{0\})$, such that $\sum_{k=1}^{\infty}\|q_k\|\leqslant+\infty$, 
there exists a unique sequence $\{p_k(t)\}$, $p_k\in{\cal O}_{\nu_k}(\overline{\cal D}\setminus\{0\})$, such that
\begin{equation}\label{eq8}
\sum_{j=0}^n a_j(t)(k+N+{\rm i}\gamma\delta_t)^j\,p_k(t)=q_k(t),
\end{equation}
and $\sum_{k=1}^{\infty}\|(k+N+{\rm i}\gamma\delta_t)^n\,p_k(t)\|\leqslant+\infty$.   


Each linear differential operator $L_k=\sum_{j=0}^n a_j(t)(k+N+{\rm i}\gamma\delta_t)^j$ is locally factorized into 
a product of the first order linear differential operators near its Fuchsian singular point $t=0$,
$$
L_k=({\rm i}\gamma)^n(\delta_t+b_1(t))\ldots(\delta_t+b_n(t)), \quad b_i\in{\cal O}(\overline{\cal D})
$$
(see, for example, \cite[\S19.3]{IY}). Furthermore, the exponents of $L_k$ at this point, $\lambda_1=-b_1(0),\ldots,\lambda_n=-b_n(0)$, are the roots of the polynomial $P(\lambda)=
\sum_{j=0}^n a_j(0)(k+N+{\rm i}\gamma\lambda)^j$ and thus are not integers (for the simplicity of notation, 
we do not put here the index $k$ to the objects related to $L_k$).

The equation (\ref{eq8}) takes the form
$$
({\rm i}\gamma)^n(\delta_t+b_1(t))\ldots(\delta_t+b_n(t))\,p_k=q_k.
$$
To understand an idea, let us study it in the case $n=1$. A general solution of the equation
$$
\delta_ty+b_1(t)y=0
$$
has the form
$$
y(t)=c\,e^{-\int\frac{b_1(t)}t\,dt}=c\,t^{\lambda_1}\tilde b_1(t), \qquad c\in{\mathbb C}, \quad 
\tilde b_1\in{\cal O}(\overline{\cal D}), \quad \tilde b_1(0)\ne0, 
$$
therefore
$$
(\delta_t+b_1(t))\,p_k=q_k
$$
implies a general expression for $p_k$:
$$
p_k(t)=c\,t^{\lambda_1}\,\tilde b_1(t)+t^{\lambda_1}\,\tilde b_1(t)
\int\frac{q_k(t)}{\tilde b_1(t)}\,t^{-\lambda_1-1}dt.
$$
Since the power exponent $\lambda_1$ is non-integer, the above expression represents a function
meromorphic at $t=0$ only if $c=0$. Expanding $q_k(t)/\tilde b_1(t)$ into a Laurent series near the point $t=0$, 
$$
\frac{q_k(t)}{\tilde b_1(t)}=\sum_{\ell\geqslant-\nu_k}q_{k\ell}\,t^{\ell},
$$
we have
$$
p_k(t)=\tilde b_1(t)\,\sum_{\ell\geqslant-\nu_k}\frac{q_{k\ell}}{\ell-\lambda_1}\,t^{\ell}\in
{\cal O}_{\nu_k}(\overline{\cal D}\setminus\{0\}).
$$
Using a similar reasoning, one concludes by induction for a general $n$ that the equation (\ref{eq8}) has a 
unique solution $p_k\in{\cal O}_{\nu_k}(\overline{\cal D}\setminus\{0\})$. To prove the convergence of 
$\sum_{k=1}^{\infty}\|(k+N+{\rm i}\gamma\delta_t)^n\,p_k(t)\|$, we write
$$
(k+N+{\rm i}\gamma\delta_t)^n\,p_k(t)=q_k(t)-\sum_{j=0}^{n-1}a_j(t)(k+N+{\rm i}\gamma\delta_t)^j\,p_k(t).
$$    
From the proof of Lemma 6 (i) it follows that
$$
\|(k+N+{\rm i}\gamma\delta_t)^j\,p_k\|\leqslant\frac1{k^{n-j}}\|(k+N+{\rm i}\gamma\delta_t)^n\,p_k\|,
$$
hence
$$
\|(k+N+{\rm i}\gamma\delta_t)^n\,p_k\|\leqslant\|q_k\|+C\Bigl(\frac1k+\ldots+\frac1{k^n}\Bigr) 
\|(k+N+{\rm i}\gamma\delta_t)^n\,p_k\|.
$$
Therefore,
$$
\|(k+N+{\rm i}\gamma\delta_t)^n\,p_k\|\leqslant\Bigl(1-\frac Ck-\ldots-\frac C{k^n}\Bigr)^{-1}\|q_k\|
$$
for $k$ large enough, and the convergence of the series $\sum_{k=1}^{\infty}\|q_k\|$ implies that of $\sum_{k=1}^{\infty}\|(k+N+{\rm i}\gamma\delta_t)^n\,p_k\|$. This proves the bijectivity of 
$\partial h/\partial\theta(0,0)$.

From the above it follows that if $\theta_1=\sum_{k=1}^{\infty}p_k(x^{{\rm i}\gamma})\,x^k\in H^n$, 
$\theta_2=\sum_{k=1}^{\infty}q_k(x^{{\rm i}\gamma})\,x^k\in H^0$ and $\partial h/\partial\theta(0,0)\theta_1=\theta_2$,
then $\|(k+N+{\rm i}\gamma\delta_t)^n\,p_k\|\leqslant A\|q_k\|$, that is, $\|\theta_1\|_n\leqslant A\|\theta_2\|_0$.
This, with Lemma 6, implies that $\partial h/\partial\theta(0,0): H^n\rightarrow H^0$ is a linear homeomorphism.

Thus, basing on the implicit mapping theorem for the mapping $h$ we conclude that there exist $\lambda_0>0$ and 
$$
\theta_{\lambda_0}=\sum_{k=1}^{\infty}p_k(x^{{\rm i}\gamma})\,x^k\in H^n
$$
such that
\begin{equation}\label{lambda_eq}
\sum_{j=0}^n a_j(x^{{\rm i}\gamma})(\delta+N)^j\theta_{\lambda_0}=\lambda_0 x\,M(\lambda_0 x,x^{{\rm i}\gamma}, \theta_{\lambda_0},\delta\theta_{\lambda_0},\dots,\delta^n\theta_{\lambda_0}).
\end{equation}
Finally, define an operator
$$
i_{\lambda_0}: \theta_{\lambda_0}\mapsto \sum_{k=1}^{\infty}p_k(x^{{\rm i}\gamma})\,(x/\lambda_0)^k
$$
from $H^0$ to ${\cal O}(S)$, $S=\lambda_0S_1$, which clearly commutes with $\delta$: 
$\delta(i_{\lambda_0}\theta)=i_{\lambda_0}(\delta\theta)$. Then applying it to the both sides of the equality 
(\ref{lambda_eq}) we obtain that the exotic series $\theta_0=i_{\lambda_0}(\theta_{\lambda_0})\in{\cal O}(S)$ is 
a solution of (\ref{eq4}). {\hfill $\Box$}

\section{An example: the third Painlev\'e equation}

Let us consider the third Painlev\'e equation
\begin{equation}\label{exeq1}
-x^2yy''+x^2(y')^2-xyy'+a\,xy^3+b\,xy+c\,x^2y^4+d\,x^2=0,
\end{equation}
where $a,b,c,d$ are complex parameters. Using general methods exposed in \cite{Br} one can find an exotic formal 
series solution $y=\varphi$ of the form
\begin{equation*}
\varphi=-\frac{4C\gamma^2 x^{{\rm i}\gamma}}{4(c\gamma^2+a^2) x^{2{\rm i}\gamma}-4aC\gamma x^{{\rm i}\gamma}+C^2}\,x^{-1}+\sum\limits_{k=1}^\infty c_k(x^{{\rm i}\gamma})\,x^k,
\end{equation*}
where $C\in\mathbb C$ and $\gamma\in\mathbb R$ are arbitrary nonzero constants and functions $c_k(t)$ are meromorphic at
$t=0$ (since the equation (\ref{exeq1}) is not only analytical but even polynomial, for applying Theorem 1 there is no requirement for an exotic formal series solution to begin with a strictly positive power of $x$).

Let us rewrite the equation \eqref{exeq1} by means of the operator $\delta$ and get the equation
\begin{equation*}
-y\,\delta^2y+(\delta y)^2+a\,x y^3+b\,x y+c\,x^2 y^4+d\,x^2=0,
\end{equation*}
that is, $F(x,y,\delta y,\delta^2y)=0$, where
$$
F(x,y_0,y_1,y_2)=-y_0y_2+y_1^2+a\,xy_0^3+b\,xy_0+c\,x^2y_0^4+d\,x^2.
$$

The partial derivatives of $F$ along the formal solution $\varphi$ are
\begin{eqnarray*}
\frac{\partial F}{\partial y_2}(x,\Phi)&=&\Bigl(\frac{4\gamma^2}C\,x^{{\rm i}\gamma}+\dots\Bigr)x^{-1}+\dots,\\
\frac{\partial F}{\partial y_1}(x,\Phi)&=&\bigl((8-8{\rm i})\,x^{{\rm i}\gamma}+\dots\bigr)x^{-1}+\dots, \\
\frac{\partial F}{\partial y_0}(x,\Phi)&=&(-8{\rm i}\,x^{{\rm i}\gamma}+\dots)x^{-1}+\dots.
\end{eqnarray*}
Since they all begin with the same power $x^{-1}$ and the order of their leading coefficient at $t=0$ (with respect to 
$t=x^{{\rm i}\gamma}$) is the same ($=1$), by Theorem 1 the series $\varphi$ converges in a sectorial domain near the origin.


\begin{thebibliography}{99}

\bibitem{Bre}
Brezhnev, Yu.\,V., A $\tau$-function solution of the sixth Painlev\'e transcendent, {\it Theor. Math. Phys.} {\bf 161}:3 
(2009), 1616--1633.

\bibitem{Br}
Bruno, A.\,D., Exotic expansions of solutions to an ordinary differential equation, {\it Doklady Math.} 
{\bf 76}:2 (2007), 729--733.

\bibitem{BrGo}
Bruno, A.\,D., Goryuchkina, I.\,V., All asymptotic expansions of solutions to the sixth Painlev\'e equation,
{\it Doklady Math.} {\bf 76}:3 (2007), 851--855.

\bibitem{BrPa}
Bruno, A.\,D., Parusnikova, A.\,V., Local expansions of solutions to the fifth Painlev\'e equation, 
{\it Doklady Math.} {\bf 83}:3 (2011), 348--352.

\bibitem{Die}
Dieudonn\'e, J. Foundations of Modern Analysis. Academic Press, New York, 1960.

\bibitem{GG1}
Gontsov, R.\,R., Goryuchkina, I.\,V., On the convergence of generalized power series satisfying an algebraic ODE,
{\it Asympt. Anal.} {\bf 93}:4 (2015), 311--325.

\bibitem{GG3}
Gontsov, R.\,R., Goryuchkina, I.\,V., The Maillet--Malgrange type theorem for generalized power series, 
{\it Manuscripta Math.} {\bf 156}:1 (2018), 171--185. 

\bibitem{GG4}
Gontsov, R.\,R., Goryuchkina, I.\,V.,  On the convergence of formal Dulac series satisfying an algebraic ODE,
{\it Sb. Math.} {\bf 210} (2019), to appear.

\bibitem{Gu}
Guzzetti, D., Tabulation of Painlev\'e 6 transcendents, {\it Nonlinearity} {\bf 25} (2012), 3235--3276.

\bibitem{Gu2}
Guzzetti, D., Poles distribution of PVI transcendents close to a critical point, {\it Physica D} {\bf 241} (2012),
2188--2203.

\bibitem{IY}
Ilyashenko, Yu.\,S., Yakovenko, S.\,Yu. Lectures on Analytic Differential Equations. Grad. Stud. Math. {\bf 86}, AMS,
2008.

\bibitem{Ma}
Malgrange, B., Sur le theor\`eme de Maillet, {\it Asympt. Anal.} {\bf 2} (1989), 1--4.

\bibitem{Sh}
Shimomura, Sh., Critical behaviours of the fifth Painlev\'e transcendents and the monodromy data,
{\it Kyushu J. Math.} {\bf 71} (2017), 139--185.





\end{thebibliography}
\end{document}